\documentclass[12pt]{amsart}
\newtheorem*{thm}{Theorem}

\newtheorem{lem}{Lemma}

\begin{document}

\title[Contractible open n-manifolds]{On covering translations and homeotopy 
groups of contractible open n-manifolds}

\author{Robert Myers}
\address{Department of Mathematics\\Oklahoma State University\\Stillwater, OK 74078}
\email{myersr@math.okstate.edu}

\thanks{Research at MSRI was supported in part by NSF grant DMS-9022140.}

\begin{abstract}
This paper gives a new proof of a result of 
Geoghegan and Mihalik which states that whenever a contractible 
open $n$-manifold $W$ which is not homeomorphic to 
$\mathbf{R}^n$ is a covering space of an $n$-manifold $M$ 
and either $n \geq 4$ or $n=3$ and $W$ is irreducible, then the 
group of covering translations injects into the homeotopy group of $W$.  
\end{abstract}

\maketitle

\section{Introduction}

The homeotopy group $\mathcal{H}(W)$ of a manifold $W$ is the group 
of homeomorphisms $Homeo(W)$ of $W$ modulo the normal subgroup 
consisting of those homeomorphisms which are isotopic to the identity. 
Given a subgroup $G$ of $Homeo(W)$ we denote by  
$\rho:G \rightarrow \mathcal{H}(W)$ the homomorphism which takes 
an element of $G$ to its isotopy class. The following theorem is a 
consequence of results of Geoghegan and Mihalik in \cite{Ge-Mi}. 

\begin{thm}[Geoghegan--Mihalik] Let $W$ be a contractible open 
$n$-man\-i\-fold which is not 
homeomorphic to $\mathbf{R}^n$. Assume that either $n \geq 4$ 
or $n=3$ and $W$ is irreducible. Suppose $W$ is a covering space of 
an $n$-manifold $M$ and $G\cong\pi_1(M)$ is the group of covering 
translations. 
Then $\rho:G \rightarrow \mathcal{H}(W)$ is one-to-one. \end{thm} 

The 3-dimensional case of this result was used by the author in 
\cite{My cover} to give examples of contractible open 3-manifolds 
which non-trivially cover other non-compact 3-manifolds but cannot 
cover compact 3-manifolds. Geoghegan and Mihalik mention this 
application in their paper but do not explicitly state the 
$n$-dimensional theorem  in the form given above. They prove much more 
general results 
in proper homotopy theory from which the theorem follows. A brief 
sketch of how to deduce the theorem from their results is given 
in section 2. 

The main result of this paper is a new, more geometric proof of this 
theorem which is inspired by David Wright's work in \cite{Wr}. 
The main ingredients are the Covering Isotopy Theorem \cite{Ce, Ed-Ki}, 
Wright's Orbit Lemma, and an analog of his Ratchet 
Lemma which we call the Isotopy Ratchet Lemma. This latter result is 
proven in section 3. The new proof of the theorem is given in section 4. 

\section{The theorem as corollary}

We first summarize some proper homotopy theory, referring to \cite{Ge-Mi} 
for further details. For our purposes we may assume that the spaces 
under consideration are locally finite simplicial complexes.  
Call a continuous function between spaces \textit{proper} if pre-images 
of compact sets are compact; two such maps are 
\textit{properly homotopic} if there is a homotopy between them which is 
a proper map. A \textit{base ray} for a non-compact space $X$ is a proper 
map $\omega :[0,\infty) \rightarrow X$. The space is 
\textit{strongly connected at $\infty$} if any 
two base rays are properly homotopic. An \textit{exhaustion} for $X$  
is a sequence $\{C_n\}_{n\geq 0}$ of compact subsets of $X$ such that 
$X=\cup_{n\geq 0}C_n$, $C_n \subseteq Int(C_{n+1})$, and every path 
component of $X-C_n$ has non-compact closure. Note that every connected, 
open PL manifold has an exhaustion. Given a base ray $\omega$ one may 
assume that the exhaustion has been chosen so that $\omega([n,\infty)) 
\subseteq X-C_n$ for all $n$.  Inclusion and change of basepoint along 
$\omega$ yield  
an inverse sequence of groups $\{\pi_1(X-C_n),\omega(n))\}$ whose limit 
$\pi_1^e(X,\omega)$ is called the \textit{fundamental group of $X$ at 
$\infty$ based at $\omega$.} Up to isomorphism this group is independent 
of the choice of exhaustion; if $X$ is strongly connected at $\infty$ it 
is also independent of the choice of base ray. 

A space $X$ is called \textit{$\pi_1$-trivial at $\infty$} 
or \textit{1-LC at $\infty$} if for 
every compact subset $A$ of $X$ there is a compact subset $B$ of $X$ 
such that $A \subseteq B$ and every loop in $X-B$ is null-homotopic 
in $X-A$; it is called \textit{simply connected at $\infty$} if, in 
addition, $B$ can be chosen so that $X-B$ is connected. Note that for 
a one-ended space (such as a contractible space) these conditions are 
equivalent. Moreover $X$ is simply connected 
at $\infty$ if and only if it is strongly connected at $\infty$ and  
$\pi_1^e(X)$ is trivial. 

Suppose $p:X \rightarrow X^{\prime}$ is a regular covering map with 
infinite cyclic group of covering translations $<g>$ generated by a 
homeomorphism $g$. We assume $X^{\prime}$ to be a locally finite 
simplicial complex. By Theorem 3.1 of Geoghegan and Mihalik \cite{Ge-Mi} 
if $X^{\prime}$ is non-compact and $g$ is properly homotopic to the 
identity, then $X$ is strongly 
connected at $\infty$; moreover if, in addition, $X$ is simply 
connected, then $\pi_1^e(X)$ is trivial, and so $X$ is simply connected 
at $\infty$. 

Now suppose that $W$ is a contractible open $n$-manifold, $n \geq 3$, 
which is a covering 
space of an $n$-manifold $M$. Consider a non-trivial cyclic subgroup 
$<g>$ of the group of covering translations; let $W^{\prime}$ be the 
quotient space $W/<g>$. Then 
$W^{\prime}$ is 
an $n$-manifold and so has the homotopy type of a finite dimensional, 
aspherical CW-complex \cite{Ml}. Thus $<g>$ must be infinite (see e.g. 
\cite[Lemma 9.4]{He}), and so $W^{\prime}$ has the homotopy 
type of a circle and is therefore non-compact. It follows that 
$W^{\prime}$ has 
a PL structure (in fact a smooth structure) \cite{Mo,Fr-Qu,Ki-Si}, 
and thus so does $W$; 
the elements of $<g>$ then preserve this structure. If $g$  
is isotopic to the identity, then by the above discussion $W$ is 
simply connected at $\infty$. But a contractible open $n$-manifold  
which is simply connected at $\infty$ is homeomorphic to $\mathbf{R}^n$ 
if either $n \geq 4$ \cite{Si,Fr} or $n=3$ and it is irreducible 
\cite{Ed,Wl}. 
(The last two references use the nominally stronger 
version of simple connectivity at $\infty$ which requires  
that $X-B$ be simply connected; for 3-manifolds the loop 
theorem implies that this is equivalent to our definition.)

\section{The isotopy ratchet lemma}

\begin{lem}[Isotopy Ratchet Lemma] Let $W$ be an open $n$-manifold 
and $g:W \rightarrow W$ a homeomorphism which is isotopic to the 
identity. Suppose $A$ is a compact subset of $W$. Then there is a 
compact subset $C$ of $W$ containing $A$ such that for all $i$, $j 
\in \mathbf{Z}$ a loop $\beta$ in $W-\cup_{m=-\infty}^{\infty} g^m(C)$ 
is null-homotopic in $W-g^i(A)$ if and only if it is null-homotopic 
in $W-g^j(A)$. \end{lem}

\begin{proof} Let $g_t$ be an isotopy with $g_1=g$ and $g_0$ equal to 
the identity. Let $C$ be a compact set whose interior contains  
$\cup_{t\in [0,1]} g_t(A)$. By the Covering Isotopy Theorem \cite{Ce, 
Ed-Ki} there is an isotopy $h_t$ such that $h_t(x)=g_t(x)$ for all 
$x \in A$ and $t \in [0,1]$, $h_t(x)=x$ for all $x \in W-C$ and $t 
\in [0,1]$, and $h_0$ is the identity. 

It suffices to prove the case $j=i+1$. Let $k_t=g^i \circ h_t \circ 
g^{-i}$. Then $k_1(g^i(A))=g^{i+1}(A)$, $k_t(x)=x$ for all $x \in 
W-g^i(C)$, and $k_0$ is the identity. 

Let $\beta$ be a loop in $W-\cup_{m=-\infty}^{\infty} g^m(C)$. 
Assume $\beta$ is null-homotopic in $W-g^i(A)$. Let $f:D \rightarrow 
W-g^i(A)$ be a map of the disk with $f(\partial D)=\beta$. 
Let $f^{\prime}=k_1 \circ f$. Then 
$f^{\prime}(\partial D)=k_1(f(\partial D))=k_1(\beta)=\beta$, and 
$f^{\prime}(D) \cap g^{i+1}(A)=k_1(f(D)) \cap g^{i+1}(A)= 
k_1(f(D)) \cap k_1(g^i(A))=k_1(f(D) \cap g^i(A))=
k_1(\emptyset)=\emptyset$. Thus $\beta$ is null-homotopic in 
$W-g^{i+1}(A)$. A similar argument, replacing $k_1$ by $k_1^{-1}$, 
establishes the converse. \end{proof} 

\section{The new proof of the theorem} 

Suppose $\Gamma$ is a group acting on an $n$-manifold $W$. It acts 
\textit{totally discontinuously} \cite{Fr-Sk} if for every compact 
subset $K$ of $W$ one has that $g(K) \cap K = \emptyset$ for all 
but finitely many elements $g$ of $\Gamma$. It acts \textit{without 
fixed points} if the only element of $\Gamma$ fixing a point is the 
identity.  Let $p:W \rightarrow W^{\prime}$ be 
the projection to the quotient space $W^{\prime}$ of the action. 
Then $\Gamma$ acts totally discontinuously and without 
fixed points on $W$ if and only if $p$ is a 
regular covering map with group of covering translations $\Gamma$ and 
$W^{\prime}$ is an $n$-manifold. (See \cite{Ma}.)

\begin{lem}[Orbit Lemma (Wright)] Let $W$ be a contractible open 
$n$-manifold, $n \geq 3$, and $g:W \rightarrow W$ a non-trivial 
homeomorphism such that the group $<g>$ of homeomorphisms generated 
by $g$ acts totally discontinuously and without fixed points on $W$. 
Given compact subsets $A$ and $C$ of $W$, there is a compact subset 
$B$ of $W$ containing $A$ such that every loop $\alpha$ in $W-B$ is 
homotopic in $W-A$ to a loop in $W-\cup_{m=-\infty}^{\infty} g^m(C)$. 
\end{lem} 

\begin{proof} This is Lemma 4.1 of \cite{Wr}. For a somewhat shorter 
proof of the case when $W$ is an irreducible 3-manifold see \cite{My free}. 
\end{proof} 

\begin{proof}[Proof of the theorem] As noted in section 2 it suffices to 
show that $W$ is $\pi_1$-trivial at $\infty$. Given a compact subset $A$ 
of $W$ and a non-trivial covering translation $g$ which is isotopic to 
the identity we let $C$ be as in the Isotopy Ratchet Lemma and then let 
$B$ be as in the Orbit Lemma. Given a loop $\alpha$ in $W-B$ we apply 
the Orbit Lemma to homotop it in $W-A$ to a loop $\beta$ in 
$W-\cup_{m=-\infty}^{\infty} g^m(C)$. 

Since $W$ is simply connected, $\beta$ is null-homotopic in $W$. Since 
$g$ is totally discontinuous $\beta$ is null-homotopic in $W-g^i(A)$ 
for some $i$, and hence by the Isotopy Ratchet Lemma it is null-homotopic 
in $W-A$. \end{proof}

\end{document}